\documentclass{amsart}
\usepackage{amscd, amssymb}

\newtheorem{theorem}{Theorem}[section]
\newtheorem{lemma}[theorem]{Lemma}
\newtheorem{proposition}[theorem]{Proposition}
\newtheorem{corollary}[theorem]{Corollary}
\newtheorem{definition}[theorem]{Definition}

\newcommand{\refT}[1]{Theorem~\ref{T:#1}}
\newcommand{\refC}[1]{Corollary~\ref{C:#1}}
\newcommand{\refP}[1]{Proposition~\ref{P:#1}}
\newcommand{\refD}[1]{Definition~\ref{D:#1}}

\newcommand{\refEq}[1]{Equation~\ref{Eq:#1}}

\newcommand{\ints}{{\mathbb Z}}

\newcommand{\reals}{{\mathbb R}}
\newcommand{\nats}{{\mathbb N}}
\newcommand{\wt}{\widetilde}
\newcommand{\vphi}{\varphi}
\newcommand{\proj}{{\mathbb P}}

\newcommand{\M}{MO^{{\mathbb Z}/2}}
\newcommand{\m}{\mathfrak{N}^{{\mathbb Z}/2}}

\begin{document}

\title{Real Equivariant Bordism and 
Stable Transversality Obstructions for $\ints/2$}
\author{Dev Sinha}
\address{Department of Mathematics, Brown University,
Providence, RI 02906}
\email{dps@math.brown.edu}

\maketitle

In this paper we compute homotopical equivariant bordism for the group
${\ints/2}$,
namely $\M_*$, geometric equivariant bordism $\m_*$, and their quotient as
$\m_*$-modules.  
This quotient is a module of stable transversality obstructions,
closely related to those of \cite{CW}.  In doing these computations, we use
the techniques of \cite{Si1}.  Because we are working in the real setting
only with $\ints/2$,
these techniques simplify greatly.  Note that $\m_*$ has been computed in
\cite{JCAlex} but the computation here of $\M_*$ is new and in fact 
simpler than that of $\m_*$, as with a modern viewpoint we make strong 
use of localization.

The paper is structured as indicated by section titles; we give basic
definitions,
then statements of theorems, and finally proofs.  Note that there is an
unevenness in the level of exposition between these sections.  The statements
of theorems are aimed more at experts in a first reading.  The proofs are
aimed at novices.

Thanks go to Peter Landweber for a close reading of this paper.

\section{Definitions}

Though all of the basic constructions can be
made for an arbitrary compact Lie group, we state things only for the group
$\ints/2$.  For a more thorough introduction see chapter fifteen in \cite{CCM}
by Costenoble.  Let $\tau$ denote the trivial 
one-dimensional real representation of ${\ints/2}$, and let $\sigma$ denote 
the non-trivial one in which the non-trivial element acts on $\reals$ by
multiplication by $-1$.  Let $BO^{\ints/2}(n)$ be the
Grassmannian of $n$-dimensional subspaces of 
$\mathcal{U} = \oplus_\infty (\tau \oplus \sigma)$, 
with $\ints/2$ action inherited from $\mathcal U$.  And let $\xi^{\ints/2}_n$
be the universal $\ints/2$ $n$-dimensional bundle over $BO^{\ints/2}(n)$,
and $T(\xi^{\ints/2}_n)$ be its associated Thom space.  Given a
representation $V$, let $S^V$ be the one-point compactification of $V$.

Let $\m_*$ denote the ring of bordism classes of $\ints/2$-manifolds (manifolds
with an involution), where bordism is defined in the usual way using
manifolds with boundary as in \cite{CF}.  To define homotopical bordism, 
first note that for any representation $V$ there are maps $S^V \wedge T(\xi^{\ints/2}_n \to  T(\xi^{\ints/2}_{n+|V|}$ defined by passage to Thom spaces
of the map classifying $V \times \xi^{\ints/2}_n$.
For non-negative $n$, $MO^{\ints/2}_n$ is the colimit 
$$ \varinjlim_W \left[ S^{n\tau \oplus W}, 
T(\xi^{\ints/2}_{|W|}) \right]^{\ints/2},$$
where $[,]^{\ints/2}$ denotes the set of ${\ints/2}$-maps, which in 
this case of taking maps from spheres is an 
abelian group.  
As in the ordinary setting, the Pontryagin-Thom construction gives rise
to a map $\m_n \to \M_n$ (again see \cite{CF}).  
But we will see that this map is {\em{not}} an isomorphism.
Indeed, we extend the definition of $\M_n$ to negative degrees in the
standard way as
$$ \varinjlim_W \left[ S^{W}, T(\xi^{\ints/2}_{|W|+|n|}) \right]
^{\ints/2}.$$
We will see that these groups are non-zero for any $n$, whereas $\m_n$
are zero for negative $n$ by definition.
We may in fact define an equivariant spectrum $\M$, equipped with deloopings
for any formal difference of representations of $\ints/2$, giving rise to 
associated homology and cohomology theories where $\M_* \cong \M_*(pt.) 
\cong MO_{\ints/2}^{-*}(pt.)$.  
There are periodicity isomorphisms (see \cite{CCM}) which imply that $\M_{V-W}$
depends only on the virtual dimension of $V-W$, so we restrict our
attention to integer gradings.

The difference between geometric and homotopical bordism arises from 
the breakdown of transversality in the presence of a group action.
The most important examples of classes in $\M_*$ not coming from geometric
bordism are the Euler classes (indeed, for $\ints/2$ we will see that
these are essentially
the only examples).  The representation $\sigma$ defines a $\ints/2$ vector
bundle over a point by projection.  There are no non-zero equivariant
sections of this bundle.  The Euler class $e_\sigma \in MO_{\ints/2}^1(pt.)$
reflects the equivariant non-triviality of this bundle.  Explicitly, given
a representation $V$ define the Euler class 
$e_V$ to be class of the 
composite $S^0 \subset S^V \to T(\xi^{\ints/2}_{|V|})$ in $\M_{-|V|}$, where the
second map is defined by passing to Thom spaces the inclusion of $V$ as
a fiber of $\xi_{|V|}$.

\section{Statements of Theorems}

Euler classes play an important role at every step in this paper.
Tom Dieck \cite{tD}, refining ideas of 
Atiyah and Segal, showed that localization
by inverting these Euler classes corresponds to ``reduction to fixed sets''.
We take tom Dieck's work as a starting point, translating his results from
the complex setting.  Once we are taking geometric fixed sets, we can make 
explicit computations.  It is through explicit computations that we 
prove the following key result.

Let $\proj(V)$ be the projective space of one-dimensional subspaces
of a representation $V$, with inherited ${\ints/2}$-action.  And let 
$[\proj(V)] \in \M_{|V|-1}$ be the image of the bordism class of 
$\proj(V)$ under the Pontryagin-Thom map.  Let $R_*$ be the sub-algebra
of $\M_*$ generated by the Euler class $e_\sigma$ and 
$[\proj(n\tau \oplus \sigma)]$, as
$n$ ranges over natural numbers.  And let $S$ be the multiplicative 
set in $R_*$ generated by $e_\sigma$.  By abuse, use $S$ to denote
the same multiplicative set in $\M_*$.

\begin{theorem}\label{T:loc}
The canonical map $S^{-1}R_* \to S^{-1}\M_*$ is an isomorphism.
\end{theorem}

In other words, any class in $\M_*$ can be multiplied by some power of
$e_\sigma$ to get a class in $R_*$ modulo the kernel of the localization
map (which is in fact 
zero by the next theorem).  Hence, to understand $\M_*$ it
suffices to understand $R_*$ and division by $e_\sigma$.  Having computed
$S^{-1}\M_*$ we can also deduce from the proof of \refT{loc} that $R_*$ is
a polynomial algebra.  In fact, $R_*$ is a maximal polynomial subalgebra
of $\M_*$.  

Remarkably, divisibility by $e_\sigma$ is
closely tied to the interplay between equivariant and ordinary bordism.
Let $\alpha \colon \M_* \to \mathfrak{N}_*$ be the augmentation map, that is
the map which forgets
$\ints/2$-action.  Clearly $e_\sigma$ is in the kernel of $\alpha$.
In some sense, this Euler class accounts for the entire kernel.

\begin{theorem}\label{T:seq}
The sequence
$$ 0 \to \M_* \overset{\cdot e_\sigma}{\to} \M_* \overset{\alpha}{\to}
\mathfrak{N}_* \to 0$$
is exact.
\end{theorem}

Using the exactness of this sequence, we define an operation on
$\M_*$.  First note that the augmentation map has a splitting 
$\iota \colon \mathfrak{N}_* \to \M_*$ which we may define by taking
a geometric representative for a class, imposing a trivial
action, and passing to homotopical bordism through the Pontryagin-Thom map.  

\begin{definition}
For any $x \in \M_*$ define $\overline{x}$ to be 
$\iota \circ \alpha(x)$.  Then $x - \overline{x}$ is in the kernel of $\alpha$,
so we define $\Gamma(x)$ to be the unique class such that $e_\sigma \Gamma(x)
= x - \overline{x}$.  
\end{definition}

Directly from \refT{loc} and \refT{seq}
we may deduce the following.

\begin{theorem}\label{T:comp}
$\M_*$ is generated over $\mathfrak{N}_*$ by $e_\sigma$ and classes
$\Gamma^i([\proj(n \tau \oplus \sigma)])$ where $i$ and $n$ range over
natural numbers.  Relations are 
\begin{itemize}
\item $e_\sigma \Gamma(x) = x - \overline{x}$
\item $\Gamma(xy) = \Gamma(x) y - \overline{x} \Gamma(y),$
\end{itemize}
where $x$ and $y$ range over the classes 
$\Gamma^i([\proj(n \tau \oplus \sigma)])$
\end{theorem}

Equivariant bordism rings have been notoriously difficult to describe.
We see here that $\M_*$ does not exhibit familiar properties.  For example, 
any multiplicative generating set of $\M_*$ must have a proper subset which
is also a multiplicative generating set.

When $x$ is a geometric class, that is $x$ is in the image of some
$[M]$ under the Pontryagin-Thom map, 
there is a geometric construction of $\Gamma(x)$, a construction which dates
back to work Conner and Floyd \cite{CF}.

\begin{definition}
Let $\gamma(M) = M \times_{\ints/2} S^1$, where $S^1 = \{(x,y) \in \reals^2|
x^2 + y^2 = 1\}$ has antipodal $\ints/2$
action and $\ints/2$ acts on the quotient by the rule 
$g \cdot [m, x, y] = [m, -x, y]$, where g is the non-trivial element
of $\ints/2$ and the brackets denote taking equivalence classes.
\end{definition}

\begin{theorem}\label{T:gamma}
$\Gamma([M]) = [\gamma(M)]$.
\end{theorem}

Because we have a geometric
model for this operation, and because we may for $\ints/2$ identify
geometric equivariant bordism as a sub-ring of homotopical equivariant
bordism, we also have explicit understanding of the geometric theory.

\begin{theorem}\label{T:geomcomp}
The geometric $\ints/2$ bordism ring
$\m_*$ is the sub-ring of $\M_*$ generated over $\mathfrak{N}_*$ by the classes
$[\gamma^i(\proj(n \tau \oplus \sigma))])$
\end{theorem}

Finally, we may identify the quotient $\M_* / \m_*$, which represents
some stable transversality obstructions, as a $\m_*$-module.

\begin{theorem}\label{T:trobs}
The $\m_*$-module $\M_*/\m_*$ is generated by classes $x_k$, $k \in \nats$,
in degree $-k$,
with relations $[\gamma(M)] x_k = [M] x_{k-1} - [\overline{M}] x_{k-1}$.
The class $x_k$ is the image under the canonical map to the quotient 
of the Euler class $e_{k \sigma}$.
\end{theorem}

Hence the Euler classes are essentially
the only transversality obstructions in this 
setting.   This is not true in general, as for $\ints/3$ for example the class
$S^\rho \overset{z \mapsto z^2}{\to} S^{\rho^{\otimes 2}} 
\to T(\xi^{\ints/3}_2)$ is a 
non-trivial class in the quotient of homotopical and geometric bordism
which is not a multiple of an Euler class.

\section{Proofs}

As should be expected, our computations start with exact sequences due to
Conner and Floyd, and tom Dieck.  

\subsection{The Conner-Floyd Exact Sequence}

Our goal is to prove the following theorem, and then explicitly compute some of
terms.

\begin{theorem}[Conner-Floyd]\label{T:CFexact}
There are maps $\eta$, $\phi$ and $\delta$, to be defined below, so that
the sequence
$$ \cdots \to \mathfrak{N}_*(B\ints/2) \overset{\eta}{\to} 
\m_* \overset{\phi}{\to} \oplus_k \mathfrak{N}_{*-k}(BO(k)) 
\overset{\delta}{\to} \mathfrak{N}_{*-1}(B\ints/2) \to \cdots
$$
is exact.
\end{theorem}

Historically, Conner and Floyd developed their exact sequence geometrically.
We proceed in what is in some sense
reverse of historical order, starting with the long exact 
sequence in $\m_*$ of the based pair $(c(E{\ints/2})_+, E{\ints/2}+)$, where $E{\ints/2}$
is a contractible space on which ${\ints/2}$ acts freely (for example,
the unit sphere in $\oplus_\infty \sigma$) and
$c(E{\ints/2})_+$ denotes the cone on $E{\ints/2}$ with a disjoint basepoint added. 
To be complete, we include the definition of relative bordism.

\begin{definition}\label{D:relbord}
Let $(X,A)$ be an admissible pair of ${\ints/2}$-spaces.  A singular
${\ints/2}$-manifold with reference to $(X,A)$ is a pair $(M, f)$ of a 
${\ints/2}$-manifold with boundary $M$ and a map $f \colon M \to X$ such that 
$f(\partial M) \subseteq A$.  Two singular ${\ints/2}$-manifolds 
$(M_1, f_1)$ and $(M_2, f_2)$ are bordant when there is a singular
${\ints/2}$-manifold $(W, g)$ such that $M_1 \sqcup M_2$ is ${\ints/2}$-diffeomorphic
to a codimension zero sub-manifold of $\partial W$, $g |_{M_1} = f_1$, 
$g|_{M_2} = f_2$, and $g(\partial W - (M_1 \sqcup M_2)) \subseteq A$.
\end{definition}

\begin{definition}
Let $\m_n(X, A)$ denote the group of equivalence classes up to
bordism of singular $\ints/2$-manifolds with reference to $(X, A)$.
\end{definition}

After elementary identifications, the exact sequence associated to 
the pair $(c(E{\ints/2})_+, E{\ints/2}_+)$ reads
\begin{equation}\label{Eq:CFexact}
\cdots \to \widetilde{\m}_*(E{\ints/2}+) \overset{i_*}{\to}
    \m_* \overset{j_*}{\to}
    \m_*(c(E{\ints/2}), E{\ints/2}) \overset{\partial}{\to} 
     \widetilde{\m}_{*-1}(E{\ints/2}+) \to \cdots.
\end{equation}

We start with analysis of the term
$\m_*(c(E{\ints/2}), E{\ints/2})$ and the map $j_*$.  
Form the bordism module of ${\ints/2}$-manifolds with free boundary, 
where we define the bordism relation as in \refD{relbord}.

\begin{proposition}
The module $\m_*(c(E{\ints/2}), E{\ints/2})$ is naturally isomorphic to the 
bordism module of ${\ints/2}$-manifolds with free boundary.
\end{proposition}

\begin{proof}
Any singular manifold with reference to $(c(E{\ints/2}), E{\ints/2})$ must 
have free boundary.  Conversely, given a ${\ints/2}$-manifold with free 
boundary, there is no obstruction to constructing a reference map to 
$(c(E{\ints/2}), E{\ints/2})$.  Applying
these observations to manifolds which play the role of bordisms, we see
that these correspondences are well-defined up to bordism
and inverse to each other.
\end{proof}

Up to bordism, a ${\ints/2}$-manifold with free boundary depends only on 
its fixed-set data.

\begin{proposition}
A ${\ints/2}$-manifold with free boundary $M$ is bordant to 
any smooth neighborhood $\mathcal{N}(M^{\ints/2})$ of the fixed set of $M$ 
as ${\ints/2}$-manifolds with free boundary.
\end{proposition}

\begin{proof}
Let $W = M \times [0, 1]$, with ``straightened angles''.  Then $\partial W$
is free outside of $M \times 0$ and $\mathcal{N}(M^{\ints/2}) \times 1$, so 
$W$ is the required bordism.
\end{proof}

Thus the map $j_* \colon \m_* \to 
\m_*(c(E{\ints/2}), E{\ints/2})$ ``reduces to fixed sets'' in the sense
of sending a representative of a bordism class 
$M$ to the bordism class of smooth neighborhoods of its fixed set.
If we choose $\mathcal{N}(M^{\ints/2})$ to be a tubular neighborhood of 
$M^{\ints/2}$ then
we can use standard equivariant differential topology to identify this tubular
neighborhood with a ${\ints/2}$-vector bundle over the fixed set, where the action 
of ${\ints/2}$ is free away from zero.  

\begin{proposition}\label{P:idenFG}
$\m_*(c(E{\ints/2}), E{\ints/2}) \cong \oplus_k \mathfrak{N}_{*-k}(BO(k))$.
\end{proposition}

\begin{proof}
Use the identifications we have made so far to equate $\m_*(c(E{\ints/2}), 
E{\ints/2})$ with the bordism module of ${\ints/2}$-vector bundles over 
trivial ${\ints/2}$-spaces where the action is free away from zero. 
The fiber of a 
${\ints/2}$-vector bundle over a trivial ${\ints/2}$-space is a representation.
Because $\ints/2$ has only one 
non-trivial representation, the action on any fiber and thus the total space
is completely determined.  Hence the forgetful map from this bordism module of 
${\ints/2}$-vector bundles which are free away from the zero section
to the the bordism module of vector bundles is an
isomorphism.  The result follows from the fact that $BO(k)$ is the classifying
space for vector bundles.  Note that we must grade according to the dimension
of the total space of the bundle in question.
\end{proof}

Interpreting the term $\wt{\m}_*(E{\ints/2}+)$ is more immediate.
From the observation that any singular manifold mapping to $E{\ints/2}$ must
itself have a free ${\ints/2}$-action, we see that this module is isomorphic
to the bordism module of ${\ints/2}$-manifolds with free ${\ints/2}$-action.
The bordism module of  
${\ints/2}$-manifolds with free ${\ints/2}$-action has a non-equivariant interpretation.

\begin{proposition}\label{P:free}
The bordism module of ${\ints/2}$-manifolds with free ${\ints/2}$-action 
is isomorphic to $\mathfrak{N}_*(B{\ints/2})$.
\end{proposition}

\begin{proof}
Consider the following diagram:
$$
	\begin{CD}
		\wt{M}		@>\tilde{f}>>	E{\ints/2}	\\
		@VVV				@VVV	\\
		M		@>f>>		B{\ints/2}.
	\end{CD}
$$
Given a representative $M$ with reference map $f$ to 
$B{\ints/2}$, pull back the
principal ${\ints/2}$-bundle $E{\ints/2}$ to get $\wt{M}$, which 
is in fact a free
${\ints/2}$-manifold.  Conversely, starting with a free ${\ints/2}$-manifold 
$\wt{M}$, there is no obstruction to constructing a map $\tilde{f}$
to $E{\ints/2}$.  Pass to quotients to obtain $f\colon M \to B{\ints/2}$.

These maps are well-defined, as we apply the previous argument to the
manifolds which act as bordisms.  
The composites of these maps are clearly identity maps.
\end{proof}

\begin{corollary}\label{C:idenNBG}
The module $\wt{\m}_*(E{\ints/2}+)$ is isomorphic to 
$\mathfrak{N}_*(B{\ints/2})$.
\end{corollary}

We may new deduce \refT{CFexact}.   We now 
give more geometric definitions of the maps in this exact sequence.

\begin{definition}
Let $\eta \colon \mathfrak{N}_*(B{\ints/2}) \to \m_*$ denote the 
$\mathfrak{N}_*$-module 
homomorphism which, using the identification of \refP{free}, sends a 
free ${\ints/2}$-bordism class to the corresponding ${\ints/2}$-bordism class.
\end{definition}

\begin{proposition}\label{P:ideneta}
The homomorphism $\eta$ coincides with the homomorphism 
$i_*$ of the exact sequence 
\ref{Eq:CFexact} under the isomorphism of \refC{idenNBG}.
\end{proposition}

\begin{definition}
Let $\vphi \colon \m_* \to \oplus_k \mathfrak{N}_{*-k} (BO(k))$ be the map
which sends a class $[M]$ to the fixed set of $[M]$ with reference map
to $\sqcup BO(k)$ classifying the normal bundle to the fixed set.
\end{definition}

\begin{proposition}\label{P:idenphi}
The homomorphism $\vphi$ coincides with the homomorphism $j_*$ of 
the exact sequence
\ref{Eq:CFexact} under the isomorphism of \refP{idenFG}.
\end{proposition}

Finally, we identify the boundary map.

\begin{definition}
Let $\delta\colon \oplus_k \mathfrak{N}_{*-k}(BO(k)) \to 
\mathfrak{N}_{*-1}(B{\ints/2})$ 
be the map of $\mathfrak{N}_*$-modules
which sends $E$, a vector bundle, to 
the sphere bundle of $E$ with fiberwise $\ints/2$ action defined by
letting the non-trivial element of $\ints/2$ act by multiplication by $-1$.
\end{definition}

\begin{proposition}\label{P:idendelta}
The homomorphism $\delta$ coincides with the homomorphism 
$\partial$ of the exact sequence of
\refEq{CFexact} under the isomorphisms of 
\refP{idenFG} and \refC{idenNBG}.
\end{proposition}

The proofs of Propositions \ref{P:ideneta}, \ref{P:idenphi} 
and \ref{P:idendelta} are straightforward. 

By Thom's seminal work, we can identify $\mathfrak{N}_*(B\ints/2)$ and 
$\oplus_k \mathfrak{N}_{*-k}(BO(k))$ given the standard computations of
the mod 2 homology of $B\ints/2$ and $BO(k)$.

\begin{proposition}\label{P:freecomp}
$\mathfrak{N}_*(B\ints/2)$ is a free $\mathfrak{N}_*$-module generated by 
classes $x_i$ in degree $i$, where $i$ ranges over natural numbers.
As the bordism module of free ${\ints/2}$-manifolds, the generator $x_i$ is
represented by the $i$-sphere with antipodal $\ints/2$-action.
\end{proposition}

\begin{proof}  
The $\mod 2$ homology of $B\ints/2 = \reals\proj^\infty$ is well known to be 
$\ints/2$ in every positive dimension.  The class in dimension $i$ is 
the image of the fundamental class of $\reals\proj^i$ under inclusion.
Under the identifications of \refP{free} these classes correspond
to spheres with antipodal action.
\end{proof}

Next note that there are  classifying maps for direct sums of the associated
universal bundles $BO(k) \times BO(l) \to BO(k+l)$.  These
maps give rise to an H-space structure on $\sqcup BO(k)$, which
in turn gives rise to a multiplication on homology.  

\begin{proposition}\label{P:compfs}
As a ring, $\oplus_k \mathfrak{N}_{*-k}(BO(k))$ is a polynomial algebra over
$\mathfrak{N}_*$
on classes $b_i \in \mathfrak{N}_{i-1}(BO(1))$ represented by the tautological
line bundle over $\reals\proj^{i-1}$.  
\end{proposition}

The fact that this ring is a polynomial algebra follows from the $\mod 2$
homology computation, which is standard.  The fact that generators are 
represented by projective spaces is a straightforward Stiefel-Whitney
number computation.

\subsection{The tom Dieck Localization Sequence}

Tom Dieck realized that there was a connection between the Conner-Floyd
exact sequence and the localization methods in equivariant $K$-theory
of Atiyah and Segal.  This
connection has been fundamentally important in our work.

The following lemma provides translation between localization and
topology.
Once again, let $S$ be the multiplicative subset of $\M_*$ generated by
$e_\sigma$.  

\begin{lemma}\label{L:inveul}
As rings, 
$\widetilde{\M_*}(S^{\oplus_\infty \sigma}) \cong S^{-1}\M_*.$
\end{lemma}

\begin{proof}
Apply $\wt{\M_*}$ to the identification $S^{\oplus_\infty \sigma} =
\varinjlim S^{\oplus_n \sigma}$.  After applying the suspension
isomorphisms $\wt{\M_*}(S^{\oplus_k \sigma}) \cong
\wt{\M_{*+1}}(S^{\oplus_{k+1}\sigma})$, the maps in the resulting
directed system are multiplication by the $e_\sigma$.
\end{proof}

Consider the cofiber sequence $S(\infty \sigma)_+ \to S^0 \to S^{\infty 
\sigma}$, which is a model for
the sequence $E{\ints/2}+ \to S^0 \to \wt{E{\ints/2}}$, 
essentially our sequence of a pair from the previous section.  
By the previous lemma, after applying $\wt{\M}_*$ to this cofiber sequence,
the second map in this sequence is the canonical map
from $\M_*$ to $S^{-1}\M_*$.  We now identify the outside terms in this
sequence.

\begin{theorem}\label{T:egtild}
$\wt{\M}_*( \wt{E{\ints/2}}) \cong 
          \oplus_{k \in \ints} \mathfrak{N}_{*-k}(BO)$.
\end{theorem}

\begin{proof}
Recall the definition of $\M_*( \wt{E{\ints/2}} )$ and consider the space
of maps from  $S^V$ to $\wt{E{\ints/2}}  \wedge T(\xi^{\ints/2}_n)$, for
any representation $V$.

First we show that for any ${\ints/2}$-spaces $X$ and $Y$, 
the restriction map 
$${\text{Maps}}^{\ints/2}(X, \wt{E{\ints/2}} \wedge Y) \to 
{\text{Maps}}(X^{\ints/2}, (\wt{E{\ints/2}} \wedge Y)^{\ints/2}) = 
{\text{Maps}}(X^{\ints/2}, Y^{\ints/2})$$ 
is an equivalence.  First note that this restriction is a fibration.  
Over a given component
of ${\text{Maps}}(X^{\ints/2}, Y^{\ints/2})$ a fiber is going to be the 
space of maps
from $X$ to $\wt{E\ints/2} \wedge Y$ which are specified on
$X^{\ints/2}$.  We filter this mapping space by filtering $X$.  Because 
the maps are already specified on $X^{\ints/2}$, we need only adjoin
cells of the form $G \times D^n$, where $G$ denotes
$\ints/2$ acting on itself by left multiplication.  Hence the 
subquotients in this filtration will be spaces of equivariant maps from  
${\ints/2} \times D^n$ to $\wt{E{\ints/2}} \wedge Y$ 
whose restriction to the boundary of ${\ints/2} \times D^n$ is specified.  
Because ${\ints/2} \times D^n$ is a free $\ints/2$-space, it suffices
to consider the restriction of such a map to one copy of $D^n$.
But $\wt{E{\ints/2}} \wedge Y$ is contractible, hence so is this mapping
space.  Therefore the fibers of our
restriction map are contractible.

Applying this argument for $X = S^V$, $Y = T(\xi^{\ints/2}_n)$ we see that
our computation follows from knowledge of 
$T(\xi^{\ints/2}_n)^{\ints/2}$.  We claim
that $T(\xi^{\ints/2}_n)^{\ints/2} = \bigvee T(\xi_i) \wedge BO(n-i)$.
We show this by analysis of the fixed set of $\xi^{\ints/2}_n$.  Any 
fixed point of $\xi^{\ints/2}_n$ must lie over a fixed point of
$BO^{\ints/2}(n)$.  But the fixed set of $BO^{\ints/2}(n)$ is the classifying
space for ${\ints/2}$-vector bundles over trivial ${\ints/2}$-spaces.  A vector bundle
over a trivial ${\ints/2}$-space decomposes as a direct sum according to decomposition
of fibers according to representation type.  As there are only two 
representations types for $\ints/2$, we deduce that $(BO^{\ints/2}(n))^{\ints/2}
= \sqcup BO(i) \times BO(n-i)$.  Restricted to a component of this fixed
set $\xi^{\ints/2}_n$ will be $\xi(i) \times \xi(n-i)$ where ${\ints/2}$ 
fixes all
points in the first factor and acts by multiplication by $-1$ on fibers
in the second factor.  Hence one component of the fixed set of 
$\xi^{\ints/2}_n$ will be $\xi(i) \times BO(n-i)$.  

Passing to Thom spaces we find 
$$T(\xi^{\ints/2}_n)^{\ints/2} = \bigvee T(\xi(i)) \wedge BO(n-i)_+$$.  
Using this result along with
our first reduction we see that 
$$[S^V, S^{\infty \sigma}\wedge T(\xi^{\ints/2}_n]^{\ints/2} = 
[S^{V^{\ints/2}}, \bigvee T(\xi(i)) \wedge BO(n-i)_+].$$
The theorem follows by passing to direct limits.
\end{proof}

As $\wt{\M}_*(\wt{E{\ints/2}}) \cong S^{-1}\M_*$ we are interested in multiplicative
structure as well.  

\begin{corollary}\label{C:loccomp}
$S^{-1}\M_* \cong \mathfrak{N}_*[x_i, e, e^{-1}]$, where as elements of
$\oplus_k \mathfrak{N}_{*-k}(BO)$,
$x_i$ are the images of the generators given in \refP{compfs} under the
canonical inclusions of $BO(k)$ into $BO$, and as an element of 
$S^{-1}\M_*$, $e$ is the image of the Euler class $e_\sigma$.
\end{corollary}

\begin{proof}
The multiplication on $\M_*$ is defined by the maps 
$T(\xi^{\ints/2}_i) \wedge T(\xi^{\ints/2}_j) \to T(\xi^{\ints/2}_{i+j})$
which are the passage to Thom spaces of the map classifying the product
of universal bundles.  These maps restrict to $T(\xi^{\ints/2}_j)^{\ints/2} =
\bigvee T(\xi(i)) \wedge BO(n-i)_+$ as the standard multiplication on 
$T(\xi(i))$ factors
smashed with the classifying map for Whitney sum on $BO(n-i)_+$ factors.

Passing to the direct limit and neglecting grading, 
we are computing $\mathfrak{N}_* (\ints \times
BO)$, where $\ints \times BO$ has an $H$-space structure which is the
product of the group structure on $\ints$ and the $H$-space structure on
$BO$ arising from Whitney sum.  The computation follows from the K\"unneth
theorem, as $\mathfrak{N}_*(\ints)$ is a Laurent polynomial ring, which
by our grading conventions is generated by a class we call $e^{-1}$ in degree
$1$, and $\mathfrak{N}_*(BO)$ is a polynomial ring in classes $x_i$ where
$x_i$ is the image of the generator of $\mathfrak{N}_i(BO(1))$ under the
inclusion from $BO(1)$ to $BO$.
That $e$ is the image of $e_\sigma$ follows directly from their definitions,
chasing through the identifications of \refT{egtild}.
\end{proof}

Next we identify the term which gives the kernel and cokernel of 
the localization map.

\begin{theorem}
$\wt{M}_*(E{\ints/2}+) \cong \mathfrak{N}_*(B{\ints/2})$.
\end{theorem}

\begin{proof}
This theorem is immediate as an application of Adams' transfer.  But in the
spirit of giving elementary proofs, we argue geometrically as follows.

From the definition of $\wt{M}_*(E{\ints/2}+)$ consider a ${\ints/2}$-map from
$S^V$ to $E{\ints/2}+ \wedge T(\xi^{\ints/2}_n)$.  Because the latter space has 
a free ${\ints/2}$-action away from the basepoint, $(S^V)^{\ints/2}$ must 
map to the 
basepoint.  If we pass to the map from the quotient $S^V / (S^V)^{\ints/2}$,
we have a ${\ints/2}$-map between ${\ints/2}$-spaces which are free ${\ints/2}$-manifolds away from
their basepoints.  Because transversality is a local condition, it is easy
to verify that transversality arguments hold in the presence of free
${\ints/2}$-actions.  Given a ${\ints/2}$-map from $S^V / (S^V)^{\ints/2}$ to 
$E{\ints/2}+ \wedge T(\xi^{\ints/2}_n)$ we may homotop it locally to a map 
which is transverse regular to the zero section of $T(\xi^{\ints/2}_n)$
and pull back a sub-manifold of $S^V / (S^V)^{\ints/2}$ which must necessarily
be free.

So following classic techniques we 
identify $\wt{\M}_*(E{\ints/2}+)$ with the bordism module of free
${\ints/2}$-manifolds.  The theorem follows from \refP{free}.
\end{proof}

As one should suspect at this point, the tom Dieck localization sequence
is precisely the homotopical version of the Conner-Floyd sequence.
The following theorem is due to tom Dieck in the complex setting.

\begin{theorem}\label{T:compare}
The diagram
$$
   \begin{CD}
	\cdots @>>> \mathfrak{N}_*(B{\ints/2}) @>>> \m_* @>\vphi>> 
\oplus_{k \in \nats} \mathfrak{N}_{*-k}(BO(k)) @>\delta>> \cdots \\
	@VVV @V\text{id}VV  @VPTVV  @V{\iota}VV \\
	\cdots @>>> \mathfrak{N}_*(B{\ints/2}) @>>> \M_* @>S^{-1}>> 
\oplus_{k \in \ints} \mathfrak{N}_{*-k}(BO) @>\partial>>
		\cdots,
   \end{CD}
$$
where the first vertical map is the identity map, the second is the
Pontryagin-Thom map and the third is defined by the standard
inclusion of $BO(k)$ into $BO$,
commutes.  
\end{theorem}

\begin{proof}
The proof of this proposition is almost immediate, as Pontryagin-Thom map $PT$
is a natural transformation of equivariant homology theories and our
exact sequence of a pair from which we defined the Conner-Floyd sequence
coincides with the cofiber sequence from which we defined the tom Dieck
sequence.
That the Pontryagin-Thom map is an
isomorphism when smashed with $E{\ints/2}+$ follows from the fact that transversality
arguments carry through in the presence of a free ${\ints/2}$-action.  
That the third vertical map is the standard inclusion follows from close
analysis of the Pontryagin-Thom map in this setting.
\end{proof}

\subsection{Proofs of The Main Results}

\begin{proof}[Proof of \refT{loc}]

We show that the images of $e_\sigma$ and $[\proj(n\tau \oplus \sigma)]$
along with $e_\sigma^{-1}$ generate $S^{-1}\M_*$.  We do so using the
explicit description of $S^{-1}\M_*$ from \refC{loccomp}.

Consider the following diagram:
$$
\begin{CD}
 	\m_* @>\vphi>>  \oplus_{k \in \nats} \mathfrak{N}_{*-k}(BO(k)) \\
	@VPTVV  @V{\iota}VV \\
	\M_* @>S^{-1}>> \oplus_{k \in \ints} \mathfrak{N}_{*-k}(BO) \cong
						\mathfrak{N}_*[x_i, e, e^{-1}],
   \end{CD}
$$
which combines results of \refT{compare} and \refC{loccomp}. 
To compute the image of 
$[\proj(n\tau \oplus \sigma)]$ under localization it suffices to look at  
fixed-set data, because it is a geometric class.

There are two components of the fixed set 
$[\proj(n\tau \oplus \sigma)]^{\ints/2}$.
Using homogeneous coordinates $[y_0, \ldots, y_n]$, 
these components are defined by the conditions
$y_n = 0$ and $z_1 = \cdots = y_{n-1} = 0$.  The condition $y_n = 0$ defines
an $n-1$ dimensional projective space.  Its normal bundle is the
tautological line bundle.  The condition $y_1 = \cdots = y_{n-1} = 0$
defines an isolated fixed point which has an $n$-dimensional normal bundle.  
Using the generators named in \refC{loccomp} we have
$$ S^{-1} [\proj(n\tau \oplus \sigma)] = x_n + e^{-(n+1)}.$$

In \refC{loccomp} we also noted that the image of $e_\sigma$ under localization
was $e$.  It thus follows that the images of $e_\sigma$ and
$[\proj(n\tau \oplus \sigma)]$ under localization, along with $e_\sigma^{-1}$,
generate $S^{-1}\M_* \cong \mathfrak{N}_*[x_i, e, e^{-1}]$, which is what was
to be shown.
\end{proof}

\begin{proof}[Proof of \refT{seq}]

The exact sequence in question is a Gysin sequence.
Apply $\wt{MO_{\ints/2}}^*$ to the cofiber sequence
$G_+ \overset{i}{\to} S^0 \overset{j}{\to} S^\sigma,$ 
where the first map is projection of $G$ onto the non-basepoint of $S^0$.
The resulting long exact sequence is
$$ \cdots \to \wt{MO_{\ints/2}}^n(S^\sigma) \overset{j^*}{\to}
{MO_{\ints/2}}^n \overset{i^*}{\to} \wt{MO_{\ints/2}}^*(G_+) 
\overset{\delta}{\to} \wt{MO_{\ints/2}}^{n+1}(S^\sigma) \to \cdots.
$$
By the periodicity of $\M$, $\wt{MO_{\ints/2}}^{n}(S^\sigma) \cong
{MO_{\ints/2}}^{n-2}$.  By definition $j^*$ is multiplication by $e_\sigma$.
From the fact that $\text{Maps}^{\ints/2}[G_+, Y]$ is homeomorphic to $Y$ for
any ${\ints/2}$-space $Y$, we see that 
$\wt{MO_{\ints/2}}^*(G_+) \cong \mathfrak{N}_*$
and $i^*$ is the augmentation map $\alpha$. 

As we remarked after the statement of \refT{seq}, the augmentation map
$\alpha$ is split.  Hence our long exact sequence breaks up into short
exact sequences, and the result follows.
\end{proof}

\begin{proof}[Proof of \refT{comp}]

First we verify relations.  Then we show that the classes listed generate
$\M_*$.  Finally we show that the relations are a complete set of relations.
It is convenient to 
view $\M_*$ as a subring of $S^{-1}\M_* \cong 
\mathfrak{N}_*[x_i, e, e^{-1}]$, which we can do as the previous
theorem implies that $e_\sigma$ is not a zero divisor.  
In this way, we may verify the second family of
relations by direct computation.  The first family of relations holds by
definition.

For convenience, rename $[\proj(n\tau \oplus \sigma)]$ as $X_n$.
By \refT{loc}, any class in $\M_*$ when multiplied by some power of
$e_\sigma$ is equal to a class in $R_*$ modulo the annihilator ideal of
$e_\sigma$, which is zero.
Hence we may filter $\M_*$ exhaustively as
$$R_* = R^0_* \subset R^1_* \subset \cdots \subset \M_*,$$
where $R^i_*$ is obtained from by adjoining to $R^{i-1}_*$ all $x \in \M_*$ 
such that $x \cdot e_\sigma = y \in R^{i-1}_*$.  By \refT{seq} the
set of all such $y$
is $\text{Ker}(\alpha) \cap R^{i-1}_*$.  The kernel of the
augmentation map is clearly generated by all classes $y - \overline{y}$. 
So we may obtain $R^i_*$ from $R^{i-1}_*$ by applying $\Gamma$ to every class
in $R^{i-1}_*$.  Since $\Gamma(xy) = \Gamma(x) y - \overline{x} \Gamma(y),$ 
it suffices to apply $\Gamma$
only to primitive elements.  It follows that $e_\sigma$ and 
$\Gamma^i(X_n)$ constitute multiplicative
generators.

Finally, to show these relations are complete we identify an additive basis
of $\M_*$.  There are two types of monomials in the additive basis, those
of the form $e_\sigma^k f$, $f \in \mathfrak{N}_*[X_i]$ and those of the
form $\Gamma^k(X_j) f$, $f \in \mathfrak{N}_*[X_i | i>j]$.  We may check that
these classes are additively independent by mapping to $S^{-1}\M_*$.
Define the {\em complication} of a monomial in our basis
elements to be the sum of the number of times both $e_\sigma$ and
$\Gamma$ appear in the monomial and the sum of all $i$ where $\Gamma^i(X_k)$
appears for some $X_k$ where $k$ is not minimal among the $X_i$ which appear.
And define the complication of a sum of monomials to be the greatest of
their individual complications.  We may use our two families of relations
to decrease complication, which inductively allows us to reduce to our
additive basis whose members have zero complication.
\end{proof}

\begin{corollary}\label{C:addbas}
An additive basis for $\M_*$ is given by monomials of
the form $e_\sigma^k f$, $f \in \mathfrak{N}_*[X_i]$ and those of the
form $\Gamma^k(X_j) f$, $f \in \mathfrak{N}_*[X_i | i>j]$.
\end{corollary}

\begin{proof}[Proof of \refT{gamma}]
Once again we use the fact that the map from 
$\M_*$ to  $S^{-1}\M_* \cong \mathfrak{N}_*[x_i, e, e^{-1}]$ is a 
faithful representaion, along with direct computation.
By \refT{compare}, we may compute the image of $[\gamma(M)]$ under
localization by analyzing fixed-set data.  

Recalling the definition of $\gamma(M)$ we see two types of fixed points
$[m, x, y]$ under the ${\ints/2}$-action, those with $x = 0$ and those with $y=0$
and $gm = m$.  The first fixed set is $\alpha(M)$, with a trivial
normal bundle.  The second fixed set is the fixed set of $M$, whose normal 
bundle is the normal bundle of this fixed set in $M$ crossed with a trivial
bundle.  By the fact that multiplication by $e^{-1}$ in 
$\mathfrak{N}_*[x_i, e, e^{-1}]$ corresponds geometrically to crossing
with a trivial bundle, this fixed set is the fixed set of $\Gamma([M])$.
\end{proof}

\begin{proof}[Proof of \refT{geomcomp}]
By analysis identical to that in the proof of \refT{seq}, the map 
$\vphi \colon \m_* \to \oplus_k \mathfrak{N}_{*-k} (BO(k))$ is injective.

We deduce from the comparison of exact sequences in \refT{compare}
that the Pontryagin-Thom map from $\m_*$ to $\M_*$ is injective.  So the
image of $\vphi$ is the image of $\M_*$ in the subring 
$\mathfrak{N}_*[x_i, e^{-1}]$ of $S^{-1}\M_*$.  The images of 
$[\gamma^i(\proj(n \tau \oplus \sigma))]$ generate this image, so
these classes generate $\m_*$.
\end{proof}

\begin{proof}[Proof of \refT{trobs}]
This theorem
follows almost immediately from \refC{addbas} and \refT{geomcomp}.
Any monomial of the form $\Gamma^k([\proj(j \tau \oplus \sigma)]) f$, 
$f \in \mathfrak{N}_*[\proj(i\tau \oplus \sigma)] | i>j]$ is in fact in 
the image of the Pontryagin-Thom map.  Monomials of the form 
$e_\sigma^k f$, $f \in \mathfrak{N}_*[\proj(i\tau \oplus \sigma)]$ are
generated over $\m_*$ by $e_\sigma^k$ which we denote by $x_k$.
The module relations for the quotient follow from the ring relations
for $\M_*$.
\end{proof}

\end{document}